%
\catcode`@=11
%
%
\def\bibn@me{R\'ef\'erences}
\def\bibliographym@rk{\centerline{{\sc\bibn@me}}
	\sectionmark\section{\ignorespaces}{\unskip\bibn@me}
	\bigbreak\bgroup
	\ifx\ninepoint\undefined\relax\else\ninepoint\fi}
%
%
%
\let\refsp@ce=\ 
\let\bibleftm@rk=[
\let\bibrightm@rk=]
%
%
%
\def\numero{n\raise.82ex\hbox{$\fam0\scriptscriptstyle o$}~\ignorespaces}
%
%
\newcount\equationc@unt
\newcount\bibc@unt
\newif\ifref@changes\ref@changesfalse
\newif\ifpageref@changes\ref@changesfalse
\newif\ifbib@changes\bib@changesfalse
\newif\ifref@undefined\ref@undefinedfalse
\newif\ifpageref@undefined\ref@undefinedfalse
\newif\ifbib@undefined\bib@undefinedfalse
\newwrite\@auxout
%
%
\def\eqnum{\global\advance\equationc@unt by 1%
\edef\lastref{\number\equationc@unt}%
\eqno{(\lastref)}}
%
%
%
%
%
%
\def\re@dreferences#1#2{{%
	\re@dreferenceslist{#1}#2,\undefined\@@}}
\def\re@dreferenceslist#1#2,#3\@@{\def\next{#2}%
	\expandafter\ifx\csname#1@@\meaning\next\endcsname\relax
	??\immediate\write16
	{Warning, #1-reference "\next" on page \the\pageno\space
	is undefined.}%
	\global\csname#1@undefinedtrue\endcsname
	\else\csname#1@@\meaning\next\endcsname\fi
	\ifx#3\undefined\relax
	\else,\refsp@ce\re@dreferenceslist{#1}#3\@@\fi}
%
%
%
\def\newlabel#1#2{{\def\next{#1}\newl@bel#2}}
\def\newl@bel#1#2{%
	\expandafter\xdef\csname ref@@\meaning\next\endcsname{#1}%
	\expandafter\xdef\csname pageref@@\meaning\next\endcsname{#2}}
\def\label#1{{%
	\toks0={#1}\message{ref(\lastref) \the\toks0,}%
	\ignorespaces\immediate\write\@auxout%
	{\noexpand\newlabel{\the\toks0}{{\lastref}{\the\pageno}}}%
	\def\next{#1}%
	\expandafter\ifx\csname ref@@\meaning\next\endcsname\lastref%
	\else\global\ref@changestrue\fi%
	\newlabel{#1}{{\lastref}{\the\pageno}}}}
\def\ref#1{\re@dreferences{ref}{#1}}
\def\pageref#1{\re@dreferences{pageref}{#1}}
%
%
\def\bibcite#1#2{{\def\next{#1}%
	\expandafter\xdef\csname bib@@\meaning\next\endcsname{#2}}}
\def\cite#1{\bibleftm@rk\re@dreferences{bib}{#1}\bibrightm@rk}
%
%
\def\beginthebibliography#1{\bibliographym@rk
	\setbox0\hbox{\bibleftm@rk#1\bibrightm@rk\enspace}
	\parindent=\wd0
	\global\bibc@unt=0
	\def\bibitem##1{\global\advance\bibc@unt by 1
		\edef\lastref{\number\bibc@unt}
		{\toks0={##1}
		\message{bib[\lastref] \the\toks0,}%
		\immediate\write\@auxout
		{\noexpand\bibcite{\the\toks0}{\lastref}}}
		\def\next{##1}%
		\expandafter\ifx
		\csname bib@@\meaning\next\endcsname\lastref
		\else\global\bib@changestrue\fi%
		\bibcite{##1}{\lastref}
		\medbreak
		\item{\hfill\bibleftm@rk\lastref\bibrightm@rk}%
		}
	}
\def\endthebibliography{\egroup\par}
%
%
%
\def\@closeaux{\closeout\@auxout
	\ifref@changes\immediate\write16
	{Warning, changes in references.}\fi
	\ifpageref@changes\immediate\write16
	{Warning, changes in page references.}\fi
	\ifbib@changes\immediate\write16
	{Warning, changes in bibliography.}\fi
	\ifref@undefined\immediate\write16
	{Warning, references undefined.}\fi
	\ifpageref@undefined\immediate\write16
	{Warning, page references undefined.}\fi
	\ifbib@undefined\immediate\write16
	{Warning, citations undefined.}\fi}
%
%
\immediate\openin\@auxout=\jobname.aux
\ifeof\@auxout \immediate\write16
  {Creating file \jobname.aux}
\immediate\closein\@auxout
\immediate\openout\@auxout=\jobname.aux
\immediate\write\@auxout {\relax}%
\immediate\closeout\@auxout
\else\immediate\closein\@auxout\fi
%
%
\input\jobname.aux
\immediate\openout\@auxout=\jobname.aux
%
%

\def\bibn@me{R\'ef\'erences bibliographiques}
%
\def\bibliographym@rk{\bgroup}
%
%
\outer\def\bye{ 	\par\vfill\supereject\end}

\def\Q{{\bf {Q}}}

\overfullrule=0pt

\magnification=1200

\def\house#1{\setbox1=\hbox{$\,#1\,$}%
\dimen1=\ht1 \advance\dimen1 by 2pt \dimen2=\dp1 \advance\dimen2 by 2pt
\setbox1=\hbox{\vrule height\dimen1 depth\dimen2\box1\vrule}%
\setbox1=\vbox{\hrule\box1}%
\advance\dimen1 by .4pt \ht1=\dimen1
\advance\dimen2 by .4pt \dp1=\dimen2 \box1\relax}

  \def\eps{{\varepsilon}}

\def\sm{\smallskip}  \def\noi{\noindent}

\def\build#1_#2^#3{\mathrel{\mathop{\kern 0pt#1}\limits_{#2}^{#3}}}

\def\date {le\ {\the\day}\ \ifcase\month\or 
janvier\or fevrier\or mars\or avril\or mai\or juin\or juillet\or
ao\^ut\or septembre\or octobre\or novembre\or 
d\'ecembre\fi\ {\oldstyle\the\year}}

\font\fivegoth=eufm5 \font\sevengoth=eufm7 \font\tengoth=eufm10

\newfam\gothfam \scriptscriptfont\gothfam=\fivegoth
\textfont\gothfam=\tengoth \scriptfont\gothfam=\sevengoth

\def\cqfd{\unskip\kern 6pt\penalty 500 \raise 0pt\hbox{\vrule\vbox 
to6pt{\hrule width 6pt \vfill\hrule}\vrule}\par}

\def\Qbar{\overline{\Q}}
\def\smallsquare{\vbox{\hrule\hbox{\vrule height 1 ex\kern 1 ex\vrule}\hrule}}
\def\cqfd{\hfill \smallsquare\vskip 3mm}

\def\mueff{{\mu_{\rm eff}}}

\def\rme{{\rm e}}


\vskip 5mm

\centerline{\bf Effective irrationality measures for real and $p$-adic roots}

\sm

\centerline{\bf of rational numbers close to $1$, with an application to}

\sm

\centerline{\bf parametric families of Thue--Mahler equations}

\vskip 13mm

\centerline{Yann B{\sevenrm UGEAUD} \footnote{}{\rm
}}

{\narrower\narrower
\vskip 15mm

\proclaim Abstract. {
We show how the theory of linear forms in two logarithms allows one to get 
effective irrationality measures for $n$-th roots of rational numbers ${a \over b}$, when 
$a$ is very close to $b$. We give a $p$-adic analogue of this result under 
the assumption that $a$ is $p$-adically very close to $b$, that is, that a large power of $p$
divides $a-b$. As an application, we solve completely 
certain families of Thue--Mahler equations. Our results illustrate, 
admittedly in a very special situation, the strength of the 
known estimates for linear forms in logarithms.
}

}

\vskip 14mm

\centerline{\bf 1. Introduction}

\vskip 6mm

Let $\xi$ be an irrational real number. The real number $\mu$ 
is an irrationality measure for $\xi$ if, for every positive $\eps$, there is
a positive number $C(\xi, \eps)$ such that every rational number ${p \over q}$
with $q \ge 1$ satisfies 
$$
\Bigl| \xi - {p \over q} \Bigr| > {C(\xi, \eps) \over q^{\mu + \eps}}.
$$
If, moreover, the constant $C(\xi, \eps)$ is effectively computable for every 
positive $\eps$, then $\mu$ is an effective irrationality measure for $\xi$. 
We denote by $\mu (\xi)$ (resp. $\mueff (\xi)$) the infimum of the 
irrationality measures (resp. effective irrationality measures) for $\xi$. 
Clearly, $\mueff (\xi)$ is larger than or equal to $\mu (\xi)$.

Every real algebraic number $\xi$ of
degree $d \ge 2$ satisfies $\mueff(\xi) \le d$, by Liouville's theorem,
and $\mu (\xi) = 2$, by Roth's theorem.
This shows that $\mueff(\xi) = 2$ if $\xi$ is quadratic, 
but the value of $\mueff (\xi)$ remains 
unknown for every $\xi$ of degree $d \ge 3$.  
Using the theory of linear forms in logarithms, Feldman \cite{Feld71}
proved that there exists a (small) 
positive real number $\tau (\xi)$ such that $\mueff(\xi) \le d - \tau (\xi)$; see 
\cite{Bu98,BiBu00} for more recent results. An 
alternative proof of Feldman's result, which does not depend on Baker's theory, was 
subsequently given by Bombieri \cite{Bo93}. 
This upper bound for $\mueff(\xi)$, valid for every real algebraic number $\xi$, 
can be considerably improved for some 
particular real algebraic numbers $\xi$, including $n$-th roots of rational 
numbers sufficiently close to $1$. This is the content of 
the following theorem of Bombieri and Mueller \cite{BoMu83}.

\proclaim Theorem BM. 
Let $a,b,n$ be positive integers with $b \ge 2$ and $n \ge 3$. 
Set 
$$
\eta :=  1 - {\log|a-b| \over \log b}.
$$ 
If $\root n \of{a/b}$ is of degree $n$ and $n > 2/ \eta$, then
$$
\mueff \Bigl( \root n \of{{a \over b}} \Bigr) \le {2 \over  \eta} +6 \Bigl(
{n^5 \log n \over \log b} \Bigr)^{1/3}.    \eqno (1.1)
$$


It follows from (1.1) that, for any positive real number $\eps$ and any $n \ge 3$, we have 
$\mueff (\root n \of{{a / b}} ) < 2 + \eps$ if $b$ is sufficiently large 
in terms of $n$ and if ${a \over b}$ is sufficiently close to~$1$.

In the same paper, Bombieri and Mueller observed that the theory of 
linear forms in logarithms implies that, for every positive real
algebraic number $\xi$, there exists
an effectively computable constant $C(\xi)$, depending only on $\xi$, such that
$$
\mueff  (\root n \of{\xi}) < C(\xi) (\log n),    \eqno (1.2)
$$
for every $n \ge 3$. When $n$ is large, (1.2) considerably improves 
Liouville's theorem. Furthermore, if $\xi$ is the rational number ${a \over b}$, 
where $a > b \ge 1$, then 
there exists an absolute, effectively computable $C$ such that 
$$
\mueff \Bigl( \root n \of{{a \over b}} \Bigr) \le C( \log a) (\log n).   \eqno (1.3) 
$$

The aim of the present note is to show how a known refinement in the theory
of linear forms in logarithms in the special
case where the rational numbers involved are very close to $1$, which goes back 
to Shorey's paper \cite{Sho74}, allows one  
to remove the dependence 
on $\log a$ in (1.3) when $a$ is very close to $b$.
Several spectacular 
applications to Diophantine problems and to Diophantine equations of 
this idea of Shorey have already been found; see 
for example \cite{Wa78,WaLiv} and the 
survey \cite{Bu08}. Quite surprisingly, it seems that it has not yet been noticed
that it can be used to give uniform, effective irrationality measures for 
roots of rational numbers and for quotients of logarithms of 
rational numbers (see \cite{Bu15}), under some suitable assumptions. 

Shorey's idea has been incorporated in the recent lower bounds for 
linear forms in Archime\-dean logarithms through a term usually denoted by $\log E$. 
Roughly speaking, the development of the 
theory of linear forms in non-Archimedean logarithms followed 
the one of its Archimedean analogue. For instance, 
the paper \cite{BuLa96} can be regarded as the $p$-adic analogue of \cite{LMN}, although
\cite{LMN} includes a parameter $\log E$ while \cite{BuLa96} does not. 
A parameter also called $\log E$ appeared for the first time 
in the $p$-adic setting in \cite{Bu99} and allows one
to get better estimates when the rational numbers involved in the 
linear form are $p$-adically close to $1$. Some applications of these refined estimates
have been given in \cite{Bu99}, a more spectacular one can be found in \cite{BeBuMi13}.
Here, we apply it to get explicit uniform, effective irrationality measures for $p$-adic
$n$-th roots of certain rational numbers.

Mignotte \cite{Mig96} was the first to observe that the introduction of the parameter 
$\log E$ in the estimates of linear forms in logarithms has a striking 
application to parametric families of Thue equations $a x^n - b y^n = c$, when 
$a$ and $b$ are positive integers very close to each other. 
A precise statement is given in Section 4. 
We extend Mignotte's result 
and solve completely multi-parametric 
families of Thue--Mahler equations. 

As will be clear in the proofs, the main results of the present note are nearly 
immediate consequences of known lower bounds for linear forms in logarithms and no new idea
is added. However, we believe that the results are striking enough to deserve 
to be pointed out. They show, admittedly in a very special situation, the strength of 
these estimates.

\vskip 5mm

\goodbreak

\centerline{\bf 2. Effective irrationality measures for real roots of rational numbers}

\vskip 6mm

Our first result gives effective irrationality measures for $n$-th roots of rational 
numbers sufficiently close to $1$. 

\proclaim Theorem 2.1.
Let $a, b, n$ be integers with $n \ge 3$ and $16 < b < a <  {6b \over 5}$. 
Define $\eta$ in $(0, 1]$ by $a - b = a^{1 - \eta}$. 
Then, we have
$$
\mueff \Bigl( \root n \of{{a \over b}} \Bigr)  
 \le {35.1 \over \eta} \, \max \Bigl\{{\log 2n  \over \eta \log a}, 10  \Bigr\}^2.
\eqno (2.1)
$$

In view of Liouville's theorem, Theorem 2.1 gives nothing new for small values of $n$ and is
only interesting for $n \ge 3511$. 

The reader may wonder whether the dependence on $n$ in (2.1) involves
$(\log n)^2$ and not only $\log n$, as in (1.3). The reason for this is that our proof
uses estimates from \cite{LMN}. We could have applied Gouillon's lower bounds \cite{Gou06}
and would have then obtained a dependence in $\log n$, but 
with much larger numerical constants; see (5.3) 
at the end of the proof of Theorem 2.1. 

We cannot deduce from Theorem 2.1 
that for any positive real number $\eps$ and any $n \ge 3$, there exist positive 
integers $a, b$ such that $\root n \of{{a \over b}}$ is of degree $n$ and has an effective 
irrationality measure less than $2 + \eps$. So, in this respect, our result is much less
interesting than Theorem BM. However, we stress that Theorem BM gives a 
non-trivial bound only when $b$ is very large compared to $n$; namely, one requires that $b$
satisfies 
$$
b > n^{216 n^2}.   \eqno (2.2)
$$
Theorem 2.1 is much stronger for smaller values of $b$. 

We point out an immediate consequence of Theorem 2.1 in the particular case $\eta = {1 \over 2}$. 

\proclaim Corollary 2.2.
Let $a, b, n$ be integers with $n \ge 3$ and $30 < b < a < b + \sqrt{a}$. 
If
$$
a \ge (2n)^{1/5}   \eqno (2.3)
$$ 
then we have
$$
\mueff  \Bigl( \root n \of{{a \over b}} \Bigr) \le 7020.
$$

The assumption (2.3) is fulfilled if $b > (2n)^{1/5}$, which is a
considerably weaker condition than (2.2).

\vskip 5mm

\goodbreak

\centerline{\bf 3. Effective irrationality measures for $p$-adic roots of rational numbers}

\vskip 6mm

Let $p$ be a prime number and $| \cdot |_p$ denote the absolute value on $\Q_p$
normalized such that $|p|_p = p^{-1}$. 
Let $\xi$ be an irrational element of $\Q_p$.
The real number $\mu$ 
is an irrationality measure for $\xi$ if, for every positive $\eps$, there is
a positive number $C(\xi, \eps)$ such that every rational number $p/q$
with $q \ge 1$ satisfies 
$$
\Bigl| \xi - {p \over q} \Bigr|_p >  {C(\xi, \eps) \over q^{\mu + \eps}}.
$$
If, moreover, the constant $C(\xi, \eps)$ is effectively computable for every 
positive $\eps$, then $\mu$ is an effective irrationality measure for $\xi$. 
We denote by $\mu (\xi)$ (resp. $\mueff (\xi)$) the infimum of the 
irrationality measures (resp. effective irrationality measures) for $\xi$. 

As in the real case, every $p$-adic algebraic number $\xi$ of
degree $d \ge 2$ satisfies $\mueff(\xi) \le d$, by Liouville's theorem,
and $\mu (\xi) = 2$, by Ridout's theorem.
This shows that $\mueff(\xi) = 2$ if $\xi$ is quadratic, 
but the value of $\mueff (\xi)$ remains 
unknown for every $\xi$ of degree $d \ge 3$.

Let $a, b, n$ be integers with $a > b \ge 1$ and $n \ge 3$. 
Let $p$ be a prime number. 
The theory of 
linear forms in $p$-adic logarithms implies that
there exists an absolute, effectively computable $C$ such that 
every $n$-th root $\zeta$ of ${a \over b}$ in $\Q_p$ satisfies
$$
\mueff (\zeta)  \le C p \log a  \log n.   \eqno (3.1) 
$$
The factor $p$ in (3.1) can be removed if 
$p$ divides $a - b$ and $|\zeta - 1|_p < 1$. 
If $p$ divides $a-b$ but does not divide $abn$, then it follows from 
Hensel's lemma that the polynomial $b X^n - a$ has a root $\zeta$ in $\Q_p$ 
such that $|\zeta - 1|_p < 1$. In the sequel, we denote this root
by $\root n \of{{a \over b}}$.


Our results concern $n$-th roots of rational numbers ${a \over b}$ 
which are $p$-adically close to $1$,
that is, such that a large power of $p$ divides $a-b$. 
We considerably improve the ($p$-adic) Liouville inequality for a class of algebraic numbers. 

\proclaim Theorem 3.1. 
Let $p$ be a prime number. 
Let $a, b$ be integers with $1 \le b < a $ and assume that $p$ divides 
$a-b$ but does not divide $ab$. 
Define $\eta$ in $(0, 1)$ by $|a - b|_p^{-1} = a^{\eta}$. 
Assume that $a^{\eta} \ge 4$. 
For any integer $n \ge 3$ which is not divisible by $p$, we have
$$
\mueff \Bigl( \root n \of{{a \over b}} \Bigr)  < 
{53.8  \over \eta} \,
\max\Bigl\{ {\log 2 n  \over \eta \log a}, 4 \Bigr\}^2
\eqno (3.2)
$$

As in Theorem 2.1, the dependence on $n$ occurs in (3.2) through the factor $(\log n)^2$. 
It is theoretically possible to reduce it to $\log n$.

We highlight an immediate consequence of Theorem 3.1 in the case 
$b=1$ and $\eta = 1/2$. 

\proclaim Corollary 3.2.
Let $p$ be a prime number.
Let $c,k, n$ be positive integers with $1 \le c < p^k$. If 
$$
p^k > (2n)^{1/4}
$$
and $p$ does not divide $n$, then we have
$$
\mueff  (\root n \of{1 + c p^k}) < 1722.
$$

\vskip 5mm

\goodbreak

\centerline{\bf 4. Parametric families of Thue--Mahler equations}

\vskip 6mm

Mignotte \cite{Mig96} was the first to observe that the introduction of the parameter 
$\log E$ in the estimates of linear forms in logarithms has a striking 
application to parametric families of Thue equations. He established, among others, the
following result.

\proclaim Theorem M. 
Let $b$ and $n$ be positive integers.
If $n$ exceeds $600$, then the only solution
in positive integers of the Thue equation
$$
(b+1) x^n - b y^n = 1  
$$
is given by $x = y = 1$.

Theorem M was subsequently improved and extended 
by Bennett and de Weger \cite{BeWe98}
in 1998. Three years later appeared a remarkable paper of Bennett \cite{Be01},
who managed to solve completely the remaining few hundreds
of Thue equations left over in \cite{BeWe98}. 

\proclaim Theorem Be. Let $a$, $b$ and $n$ be integers with
$a > b \ge 1$ and $n \ge 3$. Then, the equation
$$
| a x^n - b y^n | = 1
$$
has at most one solution in positive integers $x$ and $y$.

Theorems M and Be are closely related to Theorem 2.1, since there is a 
connection between effective irrationality measures for a given algebraic number $\xi$
and effective upper bounds for the solutions of the Thue 
equation $F(X, Y) = 1$, where $F(X, 1)$ denotes the minimal defining polynomial 
of $\xi$ over the rational integers. 

By means of  (the proof of) Theorem 3.1 we can go a step forward and solve 
parametric families of Thue--Mahler equations.

\proclaim Theorem 4.1. 
Let $s$ be a positive integer and $p_1, \ldots , p_s$ be distinct prime numbers. 
Let $\eta$ be a real number in $(0, {1 \over s+1} )$. 
Let $b \ge 2$ and $c \ge 1$ be integers such that 
$$
{\log c \over \log b} < 1 - \eta   \quad \hbox{and} \quad  
{\log |c|_{p_j}^{-1}  \over \log (b+c)} > \eta, \quad \hbox{for $j = 1, \ldots , s$}.
$$
There exists an effectively computable constant $\kappa$ such that,
for any integer $d$ with $|d| \le b$
and any integer $n$ satisfying
$$
n \ge \kappa {s \over \eta} \, \Bigl(  \log {s \over \eta} \Bigr)^2
\quad \hbox{and} \quad 
\gcd(n, p_1 \cdots p_s (p_1-1) \cdots (p_s - 1)) = 1,
$$ 
all the solutions to the Thue--Mahler equation
$$
(b + c) x^n - b y^n = d p_1^{z_1} \cdots p_s^{z_s},
$$
in integers $x, y, z_1, \ldots , z_s$ 
with $\gcd(x, y) = 1$ satisfy $|xy| \le 1$.

This is apparently the first example of a complete resolution of a 
multi-parametric family of Thue--Mahler equations.

\vskip 5mm

\goodbreak

\centerline{\bf 5. Proof of Theorem 2.1}

\vskip 6mm

We reproduce Corollaire 3 of \cite{LMN} and Corollary 2.4 of \cite{Gou06}, with 
minor simplification, in the special case
where the algebraic numbers involved are rational. 

\proclaim Theorem LMNG. 
Let $a_1, a_2, b_1, b_2$ be positive integers such that $a_1/a_2$ and $b_1/b_2$ 
are multiplicatively independent and greater than $1$. 
Let $A$ and $B$ be real numbers  such that
$$
A  \ge \max\{a_1, \rme\}, \quad   B \ge   \max\{b_1, \rme\}. 
$$
Let $u$ and $v$ be positive integers and set
$$
U'= {u  \over  \log A} + {v \over \log B}.
$$
Set 
$$
E = 1 + \min \biggl\{ {\log A  \over \log (a_1/a_2)},
{\log B \over \log (b_1 / b_2)} \biggr\}, 
$$
$$
\log U_1 = \max \{\log U' + \log  E, 600 + 150 \log E\},
$$
and
$$
\log U_2 = \max \{\log U' + \log \log E + 0.47, {10 \log E} \}.
$$
Assume furthermore that $15 \le E \le \min\{A^{3/2}, B^{3/2}\}$. Then,
$$
\log \Bigl|v \log {a_1 \over a_2} - u \log {b_1 \over b_2} \Bigr| 
\ge - 8550 (\log A) (\log B)  (\log U_1) (4 + \log E) (\log E)^{-3}   \eqno (5.1) 
$$
and
$$
\log \Bigl|v \log {a_1 \over a_2} - u \log {b_1 \over b_2} \Bigr| 
\ge - 35.1 (\log A)  (\log B)  (\log U_2)^2 (\log E)^{-3}.
\eqno (5.2)
$$

In \cite{LMN}, the authors defined the parameter $E$ 
to be equal to the right hand side of (1.8).
However, it is apparent from their proof that Theorem LMNG as stated
is correct. 

The numerical constant in (5.2) is much smaller than the one in (5.1), but the
dependence on $U'$ occurs through the factor $(\log U')^2$ in (5.2), while it
only occurs through the factor $\log U'$ in (5.1).

\medskip

\goodbreak 

\noi {\it Proof of Theorem 2.1 and of Corollary 2.2.}

Let $x, y$ be coprime integers such that $10^{10n} a < y < x \le 2y$ and
$|(a/b)^{1/n} - x/y| < 1/x^2$. 
Since
$$
\Bigl| \zeta \Bigl( {a \over b} \Bigr)^{1/n} - {x \over y} \Bigr| \le 4,
$$
for every $n$-th root of unity $\zeta$, we get
$$
\eqalign{
4^n \, \Bigl| \Bigl( {a \over b} \Bigr)^{1/n} - {x \over y} \Bigr| 
& \ge 4  \, \Bigl| {a \over b}  - \Bigl( {x \over y} \Bigr)^{n} \Bigr| \cr
& \ge   \Bigl| n \log {x \over y} - \log {a \over b} \Bigr| =: \Lambda. \cr}
$$
We apply Theorem LMNG to bound $\Lambda$ from below.


Recall that $\eta$ is defined by $a - b = a^{1 - \eta}$.
We check that
$$
{\log a \over \log (a/b)} \ge (\log a) {b \over a} a^{\eta} \ge 2.36  a^{\eta},
$$
since $a \ge 17$ and $b > 5a/6$. Since $a^{\eta} = a / (a - b) \ge 6$, we get
$2.36  a^{\eta} > 14$. Noticing that $x/y < a/b$, we get 
$$
{\log x \over \log (x/y)} \ge 2.36  a^{\eta}.
$$
Set 
$$
E := 1 + 2.36  a^{\eta}
\quad \hbox{and} \quad
U' = {n \over \log a} + {1 \over \log x}.
$$
Observe that $15 \le E \le a^{3/2}$. 
It then follows from (5.2) and the lower bound $\log E \ge \eta \log a$ that
$$
\eqalign{
\log \Lambda & \ge  - 35.1 (\log x) (\log a)  (\log E)^{-1}   
\Bigl( \max \Bigl\{{\log U' + \log \log E + 0.47 \over \log E}, 10 \Bigr\} \Bigr)^2. \cr
& \ge - {35.1 \over \eta} \, \log x \, \max \Bigl\{{\log 2n  \over \eta \log a}, 10 \Bigr\}^2, \cr}
$$
since $x$ is assumed to be sufficiently large. 

We conclude that
$$
\mueff \Bigl( \root n \of{{a \over b}} \Bigr)    
\le {35.1 \over \eta} \, \max \Bigl\{{\log 2n  \over \eta \log a}, 10  \Bigr\}^2.
$$
In particular, if 
$$
a \ge (2n)^{1/ (10 \eta)},
$$
then 
$$
\mueff \Bigl( \root n \of{{a \over b}} \Bigr)    \le {3510 \over \eta}.
$$
Choosing $\eta = {1 \over 2}$, this gives Corollary 2.2. 

Using (5.1), we obtain a better dependence on $n$, namely we get
$$
\mueff \Bigl( \root n \of{{a \over b}} \Bigr)   \le {21180 \over \eta} \,  
 \max \Bigl\{{\log (2n /  \log a) \over \eta \log a} + 1,  372 \Bigr\},    \eqno (5.3) 
$$
which is linear in $\log n$. 

\vskip 5mm

\goodbreak

\centerline{\bf 6. Proof of Theorem 3.1}

\vskip 6mm

Let $p$ be a prime number. 
Let $x_1 / y_1$ and $x_2/ y_2$ be non-zero rational numbers and
assume that there exists a real number $E$ such that 
$v_p\bigl((x_1/y_1) -1)\bigr) \ge E > 1 /(p-1)$. Theorem Bu below, established
in \cite{Bu99}, gives an explicit upper bound
for the $p$-adic valuation of
$$
\Lambda = \Bigl( \, {x_1 \over y_1} \, \Bigr)^b -
\Bigl( \, {x_2 \over y_2} \, \Bigr),
$$
where $b$ is a positive integer not divisible by $p$. 
Let $A_1 >1,A_2>1$ be real numbers such that
$$
\log A_i \ge \max \{ \log|x_i|, \log|y_i|, {E \,\log p}\},
\,\,(i=1,2).
$$
and put
$$
b' = {b \over \log A_2 } + {1 \over \log A_1}.
$$
\sm

\proclaim Theorem Bu. 
With the above notation, if $x_1 / y_1$ and $x_2/ y_2$
are multiplicatively independent, then we have the upper estimate 
$$
v_p(\Lambda ) \le {53.8  \over { E^3 \, (\log p)^4}} 
\bigl(\max \{\log b' +\log( E \, \log p) + 0.4,  4 \, E \,\log p, 5\}\bigr)^2\log A_1 \log A_2,
$$
if $p$ is odd or if $p=2$ and $v_2(x_2/y_2 - 1) \ge 2$. 

\medskip

\goodbreak 

\noi {\it Proof of Theorem 3.1 and of Corollary 3.2.}

Recall that, since $p$ does not divide $n$, every $n$-th root of unity $\zeta \not=1 $ 
in $\Qbar_p$ satisfies $v_p (\zeta - 1) = 0$. 
Let $x/y$ be a rational number. We wish to bound from above the quantity
$v_p (\root n\of {a/b} - x/y)$. Since $v_p (\root n\of {a/b} - 1)$ is positive, 
we may assume that $v_p ((x/y) - 1)$ is positive. From 
$$
\gcd \Bigl( x-y, {x^n - y^n \over x - y} \Bigr) = \gcd (x-y, n) 
$$
and the fact that $p$ does not divide $n$, we deduce that 
$$
v_p ( (x/y) - 1) = v_p ( (x/y)^n - 1).
$$
We apply Theorem Bu to bound from above the 
$p$-adic valuation of the quantity
$$
\Lambda_p :=   {a \over b} - \Bigl( {x \over y} \Bigr)^n.
$$
We introduce the parameter $E$ equal to the largest power of $p$ which divides $a-b$.
By assumption, we have $E \ge 1$ and we get 
$$
v_p ( (x/y) - 1) = v_p ( (x/y)^n - 1) = v_p ((a/b) - 1) = E.
$$
Note that $E \ge 2$ if $p=2$. 
By definition of $\eta$, we have
$$
\eta \log a = E \log p.
$$
We take $x > 10^{10n} a$ and apply Theorem Bu with
$$
\log A_1 = \max\{\log a, \eta \log a\} = \log a, \quad
\log A_2 = \max\{\log x, \eta \log a\} = \log x,
$$
and
$$
b' = {1 \over \log x} + {n \over \log a}.
$$
We then get 
$$
v_p( \Lambda_p) \le
{53.8 (\log a) (\log x) \over \eta (\log a) (\log p)} \,
\max\Bigl\{ {\log b' + \log (\eta \log a) + 0.4 \over \eta \log a}, 4, {5 \over \eta \log a} \Bigr\}^2.
$$
Since $a^{\eta} \ge \rme^{5/4}$, by assumption, this gives
$$
v_p( \Lambda_p) \le
{53.8  \log x \over \eta  \log p} \,
\max\Bigl\{ {\log b' + \log (\eta \log a) + 0.4 \over \eta \log a}, 4 \Bigr\}^2. 
$$
Thus, we obtain 
$$
v_p( \Lambda_p) \le
{53.8  \log x \over \eta  \log p} \,
\max\Bigl\{ {\log 2 n  \over \eta \log a}, 4 \Bigr\}^2,
$$
since $x$ is sufficiently large. This gives
$$
\mueff  \Bigl(\root n \of{{a \over b}} \Bigr) \le 
{53.8    \over \eta} \,
\max\Bigl\{ {\log 2 n  \over \eta \log a}, 4 \Bigr\}^2.
$$

In particular, if $a$ satisfies
$$
a \ge (2n)^{1/(4 \eta)}, 
$$ 
then we have 
$$
 \Bigl| \root n \of{{a \over b}} - \Bigl( {x \over y} \Bigr) \Bigr|_p \ge x^{-861/\eta}, 
$$
hence,
$$
\mueff  \Bigl(\root n \of{{a \over b}} \Bigr) \le {861 \over \eta}.   \eqno (6.1) 
$$
If $b=1$ and $a = 1 + c p^k$ for integers $k \ge 1$ and $c$ 
satisfying $1 \le c < p^k$, then $|a-b|_p^{-1} = p^k > \sqrt{a}$. Corollary 3.2
then follows from (6.1) with $\eta = {1 \over 2}$.

\vskip 5mm

\goodbreak

\centerline{\bf 7. Proof of Theorem 4.1}

\vskip 6mm

Let $x, y$ be integers such that
$$
(b + c) x^n - b y^n = d p_1^{z_1} \cdots p_s^{z_s}. 
$$
Assume that $|y| = \max\{|x|, |y|\}  \ge 2$, the case $|x| = \max\{|x|, |y|\}  \ge 2$ being 
analogous. Observe that
$$
\eqalign{
|d|  & = | (b + c) x^n - b y^n | \, \prod_{j=1}^s \, | (b + c) x^n - b y^n |_{p_j}  \cr
& \ge  b |y|^n \, \Bigl| \Bigl( {b+c \over b} \Bigr) \, \Bigl( {x \over y} \Bigr)^n - 1 \Bigr| 
\, \prod_{j=1}^s \, | (b + c) x^n - b y^n |_{p_j}. \cr}  \eqno (7.1)
$$
We follow the proofs of Theorem 2.1 and 3.1 to bound from below the quantities 
$$
\Bigl| \Bigl( {b+c \over b} \Bigr) \, \Bigl( {x \over y} \Bigr)^n - 1 \Bigr| 
\quad \hbox{and} \quad 
| (b + c) x^n - b y^n |_{p_j}, \quad \hbox{for $j=1, \ldots , s$}. 
$$

Unlike in the proof of Theorem 3.1, where ${x \over y}$ was assumed to be a good rational 
approximation to $\root n \of{{a \over b}}$, we have to assume that 
$n$ is coprime to $p-1$ to guarantee that ${x \over y}$
is congruent to $1$ modulo $p$, an assumption which is crucial for applying 
Theorem Bu. This 
assumption on $n$ 
implies that the $p$-adic valuation of 
$x/y - 1$ is equal to the $p$-adic valuation of $(x/y)^n - 1$.

Below, the constants implied by $\ll, \gg$ are absolute, effectively computable and positive. 
Proceeding as in the proof of Theorem 2.1, we get that
$$
\log \Bigl| \Bigl( {b+c \over b} \Bigr) \, \Bigl( {x \over y} \Bigr)^n - 1 \Bigr| 
\gg {\log y \over 1 - \eta} \max\{ \log 2n, 10 \}^2
\gg {\log y \over \eta} \max\{ \log 2n, 10 \}^2.
$$
Likewise, proceeding as in the proof of Theorem 3.1, we get 
for $j = 1, \ldots , s$ that 
$$
v_{p_j} ( (b + c) x^n - b y^n) \ll {\log y \over \eta \log p_j} \,  \max\{ \log 2n, 4 \}^2.
$$
Since $|y| \ge 2$ and $|d| \le b$, it follows from (7.1) that
$$
n \ll {(s+1) (\log 2n)^2 \over \eta}.
$$
This completes the proof of the theorem.

\vskip 8mm

\goodbreak

\centerline{\bf References}

\vskip 5mm

\beginthebibliography{999}

\bibitem{Be01}
M. A. B{ennett},
{\it Rational approximation to algebraic number of
small height~:~The Diophantine equation} $\vert a x^n - by^n \vert =1$, 
J. reine angew. Math.  535 (2001), 1--49.

\bibitem{BeBuMi13}
M. A. Bennett, Y. Bugeaud and M. Mignotte,
{\it Perfect powers with few binary digits and related Diophantine problems},
Ann. Sc. Norm. Super. Pisa Cl. Sci. (5) 12 (2013),  941--953.

\bibitem{BeWe98}
M. A. Bennett and B. M. M. de Weger,
{\it On the Diophantine equation $|a x^n - b y^n|=1$}, 
Math. Comp. 67 (1998), 413--438.

\bibitem{BiBu00}
Yu. Bilu et Y. Bugeaud,
{\it D\'emonstration du th\'eor\`eme de
Baker-Feldman via les formes lin\'eaires en deux logarithmes}, 
J. Th. Nombres Bordeaux 12 (2000), 13--23.

\bibitem{Bo93}
E. Bombieri,
{\it Effective Diophantine approximation on ${\bf G}_m$}, 
Ann. Scuola Norm. Sup. Pisa Cl. Sci. 20 (1993), 61--89.

\bibitem{BoMu83}
E. Bombieri and J. Mueller, 
{\it On effective measures of irrationality for $\root n \of{a/b}$ and related numbers}, 
J. reine angew. Math. 342 (1983), 173--196. 

\bibitem{Bu98}
Y. Bugeaud,
{\it Bornes effectives pour les solutions des \'equations en $S$-unit\'es et 
des \'equations de Thue--Mahler}, 
J. Number Theory 71 (1998), 227--244.

\bibitem{Bu99}
Y. Bugeaud,
{\it Linear forms in $p$-adic logarithms and 
the Diophantine equation
${x^n - 1 \over x - 1} = y^q$}, Math. Proc. Cambridge
Phil. Soc. 127 (1999), 373--381.

\bibitem{Bu08} 
Y. Bugeaud, 
{\it Linear forms in the logarithms of algebraic numbers
close to 1 and applications to Diophantine equations},
Proceedings of the Number Theory conference DION 2005,
Mumbai, pp. 59--76, Narosa Publ. House, 2008.

\bibitem{Bu15} 
Y. Bugeaud, 
{\it Effective irrationality measures for quotients of logarithms of rational numbers}, 
Hardy--Ramanujan  J. 38 (2015), 45--48.

\bibitem{BuLa96}
Y. Bugeaud et M. Laurent,
{\it Minoration effective de la distance $p$-adique entre puissances de
nombres alg\'ebriques}, 
J. Number Theory 61 (1996), 311--342.

\bibitem{Feld71}
N. I. Feldman,
{\it Une am\'elioration effective de l'exposant dans
le th\'eor\`eme de Liouville} (en russe), {Izv. Akad. Nauk} {35} 
(1971), 973--990. \'Egalement : 
{Math. USSR Izv.} {5} (1971), 985--1002.

\bibitem{Gou06}
N. Gouillon,
{\it Explicit lower bounds for linear forms in two logarithms}, 
J. Th\'eor. Nombres Bordeaux 18 (2006), 125--146. 

\bibitem{LMN}
{M. Laurent, M. Mignotte et Y. Nesterenko}, 
{\it Formes lin\'{e}aires en deux logarithmes et d\'{e}terminants d'interpolation}, 
{ J. Number Theory} {55} (1995), 285--321.

\bibitem{Mig96}
M. Mignotte,
{\it A note on the equation $a x^n - b y^n = c$}, 
Acta Arith. 75 (1996), 287--295.

\bibitem{Sho74} 
T. N. Shorey,
{\it Linear forms in the logarithms of algebraic numbers
with small coefficients I},
J. Indian Math. Soc. (N. S.) 38 (1974), 271--284.

\bibitem{Wa78}
M. Waldschmidt,
{\it Transcendence measures for exponentials and logarithms},
J. Austral. Math. Soc. Ser. A 25 (1978), 445--465.

\bibitem{WaLiv}
M. Waldschmidt,
Diophantine Approximation on Linear Algebraic Groups. 
Transcendence Properties of the Exponential Function in Several Variables, 
Grundlehren Math. Wiss. 326, Springer, Berlin, 2000.

\endthebibliography

\vskip1cm

\noindent Yann Bugeaud  

\noindent Universit\'e de Strasbourg

\noindent Math\'ematiques

\noindent 7, rue Ren\'e Descartes      

\noindent 67084 STRASBOURG  (FRANCE)

\vskip2mm

\noindent {\tt bugeaud@math.unistra.fr}

\bye